\documentclass[12pt,reqno]{amsart}
\usepackage{graphicx}
\usepackage{amscd}
\usepackage{amsmath}
\usepackage{amsfonts}
\usepackage{amssymb}

\theoremstyle{plain}

\numberwithin{equation}{section}

\begin{document}
\title[Encounters]{Mathematical  Encounters
}
\author{Jorge Sotomayor}
\address{Instituto de Matem\'{a}tica e Estat\'{i}stica, Univ. de  S\~{a}o Paulo \\
Rua do Mat\~ao, 1010, CEP: 05508--090, S\~{a}o Paulo, S.P, Brazil}
\email{sotp@ime.usp.br}

\thanks{The author has the partial support of  CNPq,  PQ-SR - 307690/2016-4. .
}


\maketitle
\begin{abstract} 
This 
evocative 
essay 
 focuses on mathematical activities witnessed by the author along 1962-64 at IMPA.  
The list of research problems proposed  in September 1962  by Mauricio Peixoto at the Seminar on   the
 Qualitative Theory  of Differential Equations
 is pointed out as a landmark for the beginning of the systematic research in the Qualitative Theory of Differential Equations and 
Dynamical Systems in Brazil.  \\


\end{abstract}

\section{
\label{letters}
Letters, November 1961}

\noindent  Prof. Mauricio Peixoto 

\noindent Research Institute for Advanced Study

\noindent Baltimore, Maryland, USA\smallskip

\noindent Dear Prof. Peixoto,

My name is Jorge Sotomayor. I write following 
the advice of my teacher, Prof. Jos\'e Tola\footnote {
 https://es.wikipedia.org/wiki/Jos\'e$_{-}$Tola$_{-}$Pasquel, distinguished  
 Peruvian  mathematician and educator.}, who has mentioned
to you  my interest in pursuing advanced mathematical studies in Brazil.
He has told me of your request that I send you a report on my
mathematical education, indicating, as much as possible, the level of
knowledge that I have attained so far.

Presently I am finishing the third year at the Institute of Mathematics of  San
Marcos University, Lima. Towards December I will have finished 
Mathematical Analysis (a  two--year course, which
covers most of the book by T. Apostol). 
The same for Differential Geometry (a
year course, which covers almost completely the book of  T. Willmore) and
Analytic Functions (a  one--semester course, based on the book of Ahlfors
and notes by Prof. Tola). 
Last semester I took a course on Modern Algebra,
also taught by Tola, based on parts of  van der Waerden and local lecture notes.

I mention only the courses to which I have devoted special attention, and in 
which I feel somewhat confident and fulfilled. I particularly enjoy the books of
Apostol and Willmore and have worked a lot on them.

So far I have remained ignorant of the theoretical aspects of Differential
Equations, having only taken last year a one--semester course on  methods of
solutions of ordinary ones. No other courses on ODEs or PDEs are taught here.

\hfill J. S.

\smallskip

\noindent Dear Sotomayor,

Next year, towards April, I will start at IMPA a series  of courses and
seminars on  the Qualitative Theory of Differential Equations. In order
that you  may 
have a chance to catch up with the level of the other prospective
participants, I suggest you do the following reading:

On \emph{Differential Equations}: study \emph{Coddington and Levinson},
Chapters 1, 3, 15 and 16. You might find more accessible the little book by
\emph{Hurewicz}, whose contents almost coincide with  the material listed above.

On \emph{Topology}: read \emph{Seifert and Therfall}; you should try to reach
Homology, the structure of surfaces, and the Euler Characteristic on manifolds.
The little book by \emph{Pontryagin} is also suitable for this.

I am writing now to Dr. L\'{e}lio Gamma, Director of IMPA, recommending you
for a fellowship. You will hear directly from him later.

\hfill M. P.

\smallskip
The fellowship was granted. 
This remarkable gesture of solidarity made 
 possible my arrival to Rio de Janeiro on time for the beginning of 
the seminar  on the Qualitative Theory of Differential Equations.
 To this end the reading suggestions in the letter  above  were faithfully followed,  not  without  considerable struggle.

\smallskip
 From   December 1961  to February   1962,  took place  my mathematical preparation for the trip to Brazil.  The  outcome of this endeavor was the   
 feeling of the   existence of  a Theory of Differential Equations founded, not in only in Calculus, but in Mathematical Analysis and in Geometry. The Poincar\'e - Bendixson Theorem was a remarkable sample. 
 This initial sensation  
 led  me to  fancy that   things of great  interest could dwell  inside such Theory.  

\section{A List of  Open Problems in ODEs, September 1962.  \label{list_english}} 
 Besides
 Peixoto's
 Seminar, other engaging 
 mathematical activities   were witnessed  by the author at IMPA along  the years  1962-64. 
In 
\cite{l},  on which this section is based, 
 was given a panoramic view of them as well as of the preparatory  topics presented at the Seminar,  previous to Peixoto's lecture: 
\emph{Open Problems in the Qualitative Theory of Ordinary Differential
Equations,}    delivered  in September \, 1962.  

\smallskip

 The following problems were proposed and  discussed:
 
\smallskip
1. \emph{First order structurally stable systems}. 

Consider the complement
$\mathcal{B}$ of the set $\mathcal{S}$ of \ $C^{r}$-structurally stable
vector fields on a two-dimensional manifold, relative to the set 
$\mathcal{X}$ of all vector fields on the manifold. Let 
$\mathcal{B}$ be endowed with the $C^{r}$ topology.
Characterize the set $\mathcal{S}_{1}$ of those vector fields that are
structurally stable with respect to arbitrarily small perturbations inside 
$\mathcal{B=X\backslash S}%
$.\smallskip\ 

This problem goes back to a 1938 research announcement of A. A. Andronov and 
E. A. Leontovich 
\cite{al_38}, \cite{al}. 
They gave a characterization of $\mathcal{S}_{1}$ 
for a compact disk in the plane. This initiates the
systematic study of the bifurcations (qualitative changes)
that occur in families of vector fields as they cross $\mathcal{B}$.
In the research announcement  --contained in a dense four pages note--  they stated that the most stable
bifurcations, regarding $\mathcal B$ with the induced topology,  occur in $\mathcal{S}_{1}$,   \cite{s}, \cite{scd}. 

\smallskip

2. \emph{The problem of the arc}. 

\ Prove or disprove that a continuous curve (an arc)
in the space $\mathcal{X}$ of vector fields of class $C^{r}$ on the sphere can
be arbitrarily well-approximated by a continuous curve that meets only 
finitely many \emph{bifurcation points}; 
that is, points  outside  the set of 
structurally stable vector fields, at which qualitative changes occur.

Later research
established that both   
 share  great complexity,
 which grows quickly with   the dimension. 
This  knowledge 
became 
apparent after the work of 
S. Smale \cite{BAMS} 
and also Newhouse  \cite{n}  and Palis -Takens  \cite{p_t}.  

The  understanding of the phenomenon of persistent accumulation of bifurcations implies that  the problem of the arc as stated above
 has  a negative answer. 
However, after removing the requirement of the approximation, 
Peixoto  
and
S. Newhouse proved that every pair of structurally stable vector fields
is connected by an arc that meets only finitely many bifurcation points. See \cite{p_n}.

\smallskip

3. \emph{The classification problem}. 

Use combinatorial invariants
to classify the connected components of the open set of structurally stable 
vector fields. The essential
difficulty of this problem is to determine when two structurally stable 
vector fields agree up to a homeomorphism that preserves their orbits and is isotopic to the
identity.
\smallskip

Progress on this problem has been made by C. Gutierrez (1944 - 2008), W. de Melo (1947 -  2016), and by
Peixoto   \cite{p_s}. 

\smallskip
\smallskip

4. \emph{The existence of nontrivial minimal sets}.  

Do invariant perfect sets (that is, sets that are nonempty, compact, and 
transversally totally discontinuous) exist for differential equations of 
class $C^{2}$ on orientable two-dimensional manifolds?

\smallskip
This problem goes back to H. Poincar\'{e} and A. Denjoy and was 
known to experts. 

It was solved in the negative direction by A. J. Schwartz \cite{sch_d}. Peixoto presented this result
from a preprint 
that he received in November 1962. 

\smallskip

5. \emph{Structurally stable second order differential equations}. 

For equations of the form $x^{\prime\prime}=f(x,x^{\prime})$ 
(more precisely, for systems of the form $x^{\prime}=y,\,\,y^{\prime}=f(x,y)$), 
characterize structural stability, and prove the genericity of 
structural stability,  in the spirit  of  Peixoto's results for vector fields
on two-dimensional manifolds.
\smallskip\

Problems 2 to 4  were assigned, 
in one-to-one correspondence, to the
senior participants of  the seminar. 
The first and last  problems   were held in reserve for a few months. 

\smallskip

Peixoto's support led me to attack the first problem
on his list.

The preparation for the presentation
of the  understanding the results,  formulated   for plane regions,  by Andronov and Leontovich, \cite{ al_38}, and the subsequent struggle 
to extend them to  surfaces,   that  I made in  May 1963  at  Peixoto's seminar,  led   me to  the  
 hunch 
that  the First and Second Problems were intimately interconnected.   \\
 The effort to provide a conceptual geometric synthesis  for my extension of \cite{ al_38} to  surfaces
 was  further elaborated for  my participation in a session of short communications delivered 
  at the Fourth Brazilian Mathematics Colloquium, in   July 1963, \cite{cbm_63} , \cite{l}, where all the people  
 working  under Peixoto's supervision 
 made  reports. 
 Only  Ivan Kupka, the most knowledgeable of the group,   delivered a plenary lecture.
  
 The  critical appreciation of  the {\it ensemble}  of these presentations   led  me  to  the strong conviction  that  the First and Second Problems were part of the same problem and that their   suitable synthesis should be presented in terms of infinite dimensional submanifolds and transversality in the Banach space of all vector fields on a surface, \cite{l}.

This was the beginning of a mathematical 
endeavor
that
would naturally led to the 
second one and, later,  touched also the
fifth one \cite{s2}.

The results of the  work carried out by the author,
 inspired   by 
 the
geometrization of differential equations in 
Peixoto's papers  ~\cite{p} and ~\cite{k},  
were  presented in \cite{s}.  
\section{ A Glimpse into Structural Stability 
}

The concept of structural stability was established during the collaboration
of the 
Russian 
mathematicians A. Andronov and L. Pontryagin that
started in 1932 \cite{pont_auto}. It first appeared in their research note published in 1937. 
Andronov (who was also a physicist) 
started  
a very  important Russian  school in Dynamical Systems. 
He left a
remarkable mathematical heritage, 
highly respected both in Russia and in the West. 

For a dynamic model --that is, a differential equation or system 
$x^{\prime}=f(x)$-- to faithfully represent a phenomenon of the physical
world, it must have a certain degree of stability. Small perturbations,
unavoidable in the recording of data and experimentation, should not affect
its essential features. Mathematically this is expressed by the requirement
that the \emph{phase portrait} of the model, 
which is the geometric synthesis of the system, 
must be topologically unchanged by small perturbations. 
In other words, the phase
portraits of $f$ and $f+\Delta f$ must agree up to a homeomorphism of the
form $I+\Delta I$, where $I$ is the identity   
transformation of the phase space of the system and $||\Delta I||$ is small. 
A homeomorphism of the form $I+\Delta I$ is called an
$\epsilon$-{\emph homeomorphism}
if $||\Delta I||<\epsilon$; that is, it moves points
at most $\epsilon$ units from their original positions.

Andronov and Pontryagin stated a characterization of structurally stable
systems on a disk in the plane. 
This work was supported by the analysis of numerous
concrete models of mechanical systems and electrical circuits, performed by Andronov and his associates \cite{a}. The
concept of structural stability, initially called \emph{robustness}, 
represents a remarkable evolution of the 
continuation method of Poincar\'{e}. 

When the 
American mathematician S. Lefschetz translated the
writings of Andronov and his collaborators from Russian to English \cite{a_L}, 
he changed the name of the concept to the more
descriptive one it has today. He also stimulated 
H. B. de Baggis
to work on a proof of the main result as stated by Andronov and
Pontryagin.

Peixoto improved the results of the Russian pioneers in several directions. 

For example, he introduced the space $\mathcal{X}$ of all vector fields,
and he established the 
openness and genericity of  structurally stable vector fields on the 
plane and on orientable surfaces. He also removed the $\epsilon$-homeomorphism  
requirement from the original definition, proving that it is equivalent to the  existence of any homeomorphism. This was a substantial improvement of 
the Andronov-Pontryagin theory, which was 
formulated  for a disk in the plane \cite{p}.

The transition from the plane to surfaces, as  
in 
Peixoto's work,  
takes us from classical ODEs to  the modern theory
of Dynamical Systems,   from Andronov and Pontryagin to D. V. Anosov and S. Smale. 

It has also raised 
delicate problems --for instance, the 
\emph{closing lemma}--  that have challenged mathematicians 
for decades \cite{g}.

\smallskip

In \cite{sm} S. Smale regards Pei\-xo\-to's structural stability theorem as the 
prototypical example and fundamental model to follow for 
global analysis. 

In Sotomayor \cite{sw} the reader can find 
biographical data about Mauricio Peixoto.     

\smallskip

\section{ Other Encounters, some reaching a more distant past}

Once a general surface, not just the plane, is regarded as the natural 
domain for the global analysis of differential equations, it is natural to
inquire into the stability properties of the differential equations of 
Classical Differential Geometry that are 
naturally  attached  to 
curved surfaces.  Examples are given by the differential equations 
defining the lines of principal curvature, see  \cite{gs},  and the  asymptotic curves, see \cite{sg}. 

This shift in perspective assigned a new significance to
questions that are suggested by  attentive readings 
after fortunate  encounters with
 G. Monge, C. Dupin, 
and G. Darboux (see \cite{m}, \cite{sg}, and \cite{gs}).

Elaboration  of the
some ideas of  Classical Differential Geometry
were 
recast  in geometric language  and   in  a new perspective, under the influence of ideas such as Genericity, Structural Stability and Bifurcation Theory, coming from Differential Equations and Dynamical Systems. 
A partial view of the results obtained can be found in 
\cite{m}, \cite{sg}, and \cite{gs}.
A recent historical essay on this subject can be found in Garcia - Sotomayor \cite{A_M}.

The  encounters  with Peixoto, epistolar,  on September\ , 1961 and, in presence,  at his  lecture on the 
list of ODE problems in September \, 1962,  
were the beginning  of several other equally  fruitful ones. 

I will mention only two of them: Steve Smale (1930 - \; ), Berkeley,  1966-67 and Ren\'ee Thom (1923 - 2004), Bures sur Ivette September - December, 1972. The deep contributions of these  outstanding mathematicians are still in the process of as\-si\-mi\-la\-tion.  

I attribute to the influence of Smale my initial  incursion in n-di\-men\-sion\-al bifurcations \cite{bif_n}, about which  I made a presentation toward November  1987 in  his  inspiring and challenging   Dynamical Systems Seminar. 

To Thom I owe my initial interest in the Theory of Singularities \cite{sing_76}.

\smallskip
\section{The Beginning  of the Brazilian School in Dynamical Systems}

The  Open Problems session in September  1962 
was the high point  of the seminar that Peixoto inaugurated  at IMPA 
in April of that year. 

The evocative  essay \emph{A List of Problems on ODEs} \cite{l} 
 sets 
 Peixoto's lecture  reported above 
  in  
  the center of the stage for the mathematical research activities at 
IMPA --its courses, seminars, and visitors-- during the years 1962--64. 

Peixoto's endeavors as Research Director and his
seminar on the Qualitative Theory of Differential Equations 
launched the first effective effort in Brazil to stimulate 
research in the field.

This  landmark in the history of Brazilian mathematics constitutes the
beginning of the Brazilian school of Dynamical Systems. 

After the starting step given by Peixoto, several successive   generations  in Dynamical Systems,  in diversified research directions came in, 
spreading in Brazil and abroad.

The achievement of Arthur  Avila, recipient  of  the  Fields Medal 2014,  must be pointed out as a
remarkable  landmark on this field  \cite{A_G}.
  
\smallskip

\noindent \textbf{Acknowledgement} The author is grateful to L. F. Mello and R. A. Garcia for helpful comments.

\end{document}